\documentstyle[12pt,amssymb]{article}
\newtheorem{theorem}{Theorem}[section]
\newtheorem{proposition}[theorem]{Proposition}
\newtheorem{lemma}[theorem]{Lemma}
\newtheorem{definition}[theorem]{Definition}
\newtheorem{corollary}[theorem]{Corollary}

\newcommand{\bref}[1]{(\ref{#1})}

\newcommand{\half}{{\frac{1}{2}}}
\renewcommand{\l}{\lambda}
\newcommand{\e}{\epsilon}

\renewcommand{\P}{{\cal P}}

\newcommand{\G}{\Gamma}

\renewcommand{\d}{\delta}

\newcommand{\R}{{\Bbb R}}
\newcommand{\Z}{{\Bbb Z}}

\newcommand{\T}{{\cal T}}

\newcommand{\lo}[1]{{(\log\frac{1}{\d})^{#1}}}

\newcommand{\supp}{\hbox{supp\,}}
\newcommand{\qed}{\mbox{\hfill$\square$}}
\linethickness{.4mm}

\begin{document}

\title{A local smoothing estimate in higher dimensions}

\author{Izabella {\L}aba and Thomas Wolff}
 
\maketitle 
 
The purpose of this paper is to prove the higher-dimensional analogue of
the local smoothing estimate of \cite{Wls}.  

We denote by $\Gamma$ the forward light cone 
\[
\Gamma=\{\xi\in\R^{d+1}:\ \xi_{d+1}=\sqrt{\xi_1^2+\dots+\xi_d^2}\},
\]
Let $N$ be a large parameter, $C$ a constant, and let 
$\Gamma_N(C)$ denote the $C$-neighborhood of the
cone segment $\{\xi:\ 2^{-C}N\leq |\xi|\leq 2^CN\}$.
For fixed $N$, we take a partition of unity subordinate to a covering of
$S^{d-1}$ by caps $\Theta$ of diameter about $N^{-\frac{1}{2}}$, and use this
to form a (smooth) partition of unity $y_{\Theta}$ on $\G_N(C)$ in the natural way.
We will write $\G_{N,\Theta}(C)=\supp y_{\Theta}$.
Let $\Xi_{\Theta}$ be a function whose Fourier transform coincides with
$y_{\Theta}$ on $\G_N(1)$. If 
the support of $\hat{f}$ is contained in $\G_N(1)$, we define
\[
\|f\|_{p, mic}=\left(\sum_{\Theta}\|\Xi_{\Theta}\ast f\|_p^p\right)^{\frac{1}{p}}
\]
for $2\leq p<\infty$, and
\[
\|f\|_{\infty, mic}=\sup_{\Theta}\|\Xi_{\Theta}\ast f\|_{\infty}.
\]

\begin{theorem}\label{thm1}
The following  estimate holds if $d\geq 3$,
$p> p_d\stackrel{def}{=}\min(2+\frac{8}{d-3}, 2+\frac{32}{3d-7})$ and $\hat{f}$
is supported in $\Gamma_N(1)$:
\begin{equation}
\forall\e\,\exists C_{\e}:\;\|f\|_p\leq C_{\e}N^{\e}N^{\frac{d-1}{2}-\frac{d}{p}}
\|f\|_{p, mic}.
\label{qa3}\end{equation}
\end{theorem}

This is sharp except for endpoint issues for the indicated values of $p$, but
on the other hand the expected range of $p$ is $p\geq 2+\frac{4}{d-1}$; see
the introduction to \cite{Wls}. We opt in  this paper  for simplicity over
efficiency. As will be seen, Theorem \ref{thm1} is much easier than its two-dimensional
analogue proved in \cite{Wls}; in particular, the geometrical arguments
involved are much simpler.  Improvements in the exponent should be possible,
for example by extending the geometrical analysis of Section 1 of \cite{Wls}
(see also \cite{M}, \cite{Sc}, \cite{AmerJ}, \cite{surv}) to higher dimensions.
Since this would complicate
the paper considerably and could not in any case give a sharp result,
we decided against carrying it out here.

Theorem \ref{thm1} implies the following partial result on the $d+1$-dimensional 
local smoothing and cone multiplier problems. We let $\|f\|_{p,\alpha}$ be
the inhomogeneous Sobolev norm with $\alpha$ derivatives in $L^p$.

\begin{corollary}\label{cor2}
(i) If $u$ is a solution of $\square u=0$, $u(\cdot, 0)=f$,
$\frac{du}{dt}(\cdot, 0)=0$ in $\R^{d+1}$ then
\[
\|u\|_{L^p(\R^{d}\times[1,2])}\leq C_{p\alpha}\|f\|_{p,\alpha}
\]
if $p> p_d$ and $\alpha>\frac{d-1}{2}-\frac{d}{p}$. 

\smallskip

(ii) Let $\rho_1$ be a function in $C_0^{\infty}((1,2))$, and let
$\rho_2\in C_0^\infty(\R^{d})$. 
Then the cone multiplier operators $T_{\alpha}$ defined via 
$\widehat{T_{\alpha}f}=m_{\alpha}\hat{f}$, where 
\[
m_{\alpha}(x)=|x_{d+1}-\sqrt{x_1^2+\dots+x_d^2}|^{\alpha}\rho_1(x_{d+1})
\rho_2(x_1,\dots,x_d),
\]
are bounded on $L^p(\R^{d+1})$ if $p>p_d$ and $\alpha>\frac{d-1}{2}-\frac{d}{p}$.

\end{corollary}

The proof is identical to that of Corollary 2 in \cite{Wls} and we will not
reproduce it here.  For further discussion we refer the reader to
e.g. \cite{B}, \cite{MSS}, \cite{Stein}, \cite{TV}, \cite{bil}, \cite{Wls}.

The proof of Theorem \ref{thm1} follows the general outline of the proof
of the $2+1$-dimensional result in \cite{Wls}.  We will in particular
rely on the ``induction on scales" argument of \cite{Wls}: assuming that
the estimate \bref{qa3} is known on scale $\sqrt{N}$, we can prove it
on scale $N$ by applying it once on scale $\sqrt{N}$ and once on a
slightly smaller scale $N^{\frac{1}{2}-\e_0}$ for some $\e_0>0$.
As in \cite{Wls}, the crucial step of passing from scale $\sqrt{N}$ to 
$N^{\frac{1}{2}-\e_0}$ uses a certain {\em localization property} of
functions, proved using geometrical arguments.  The geometry and
combinatorics involved is, however, much simpler than in the $d=2$ case;
in particular, instead of the complicated bounds on circle tangencies 
proved in \cite{Wls} we only need a fairly simple lemma concerning 
incidences between a set of points and a family of separated ``plates" or tubes.
We remark that the proof of Theorem \ref{thm1} for $d\geq 4$ and 
$p_d=2+\frac{8}{d-3}$ could be simplified even further, as it requires
changing scales only once.

The paper is organized as follows.  In Section \ref{prel} we explain the
notation used throughout the paper and prove some general properties
of the norms $\|\ \|_{p,mic}$.  We then introduce {\em $N$-functions},
which are -- in a useful sense -- the basic components of functions
with Fourier support in $\G_{N}(C)$ (Section \ref{sec4}).  In Section \ref{sec5}
we deduce Theorem \ref{thm1} from the inductive step in Proposition
\ref{scales}.  The rest of the article is devoted to the proof of 
Proposition \ref{scales}.  In Section \ref{sec2} we obtain the necessary
geometric information, including the incidence lemmas mentioned above.
This information is used in Section \ref{sec3} to obtain the required
localization properties of $N$-functions.  The main induction on scales
argument is given in Section \ref{sec5a}.

This work was partially supported by NSERC grant 22R80520 and by 
NSF grant DMS-0105158.


\section{Notation and preliminaries}\label{prel}


Throughout this paper we will fix the value of $d\geq 3$.  We will use
$N$ to denote a large parameter and $\delta$ to denote a small parameter; 
unless specified otherwise, we will assume that $\d=N^{-1}$.  All constants appearing 
in the sequel, including $\e_0$, $\e$, $M$, $C_i$, will depend on $d$ and $p$
but not on $N$ or $\d$.  We will write $A\lesssim B$ if $A\leq C B$ with
the constant $C$ independent of $N$ and $\d$, and $A\approx B$ if $A\lesssim B$
and $B\lesssim A$.  We will also write $A\lessapprox B$ if $A\lesssim \lo{C}B$
for some constant $C$.
The constants $C,C_i$, and the implicit constants in
$\lesssim$ and $\lessapprox$ will be adjusted numerous times throughout the
proof, in particular after each application of Proposition \ref{scales}.
The constants $\e_0$ will be assumed to be sufficiently small
and will remain constant throughout the proof; we also let $0<\e<\e_0^2$.
Except when specified otherwise, $t$ will be a dyadic number such that
$t\approx \d^{\e_0}$.

We will use $\chi_E$ to denote the indicator function of the set $E$, and
$|E|$ to denote the Lebesgue measure or cardinality of $E$ depending on the context.
A {\em logarithmic fraction of $E$} will be a subset of $E$ with measure $\gtrapprox|E|$. 


Let a family of sets ${\cal S}_N$ be given for each $N$. We will say that
${\cal S}_N$ have {\em finite overlap} if there is a constant $C$ such
that for any $N$ any point in $\R^{d+1}$ belongs to at most $C$ sets
in ${\cal S}_N$. 

A {\em $\d$-plate} is a rectangular box of size $C_0\d\times C_0\d^{\frac{1}{2}}
\times\dots\times C_0\d^{\frac{1}{2}}\times C_0$ whose longest axis is a light ray
and whose axes of length $C_0\d^{\frac{1}{2}}$ are tangent to the corresponding
light cone.  A {\em $\d$-tube} is a rectangular box of size 
$C_0\d\times\dots\times C_0\d\times C_0$ whose longest axis is a light ray.
The {\em direction} of a tube or plate is the direction of its longest axis.
This direction will always be of the form $(e,1)$ with $e\in S^{d-1}$. 
Two $\d$-plates or $\d$-tubes are {\em comparable}
if one is contained in the dilate of the other by a fixed constant $C$, and
they are {\em parallel} if their axis directions fail to be $C\d$-separated
for a suitable  $C$. A family of $\d$-plates or $\d$-tubes is 
{\em separated} if no more than $C$ are comparable to
any given one. 

If $\pi$ is a $C_0\times C_0\d^{\frac{1}{2}}\times\dots\times C_0\d^{\frac{1}{2}}
\times C_0\d$-plate, with respective axes $e_1,\dots,e_{d+1}$, then we let 
$\pi^{*}$ be a rectangle centered at the point $Ne_{d+1}$ with axes
$e_1,\dots,e_{d+1}$ and respective axis lengths $C_1, C_1N^{\frac{1}{2}},
\dots,C_1N^{\half}, C_1N$, where $C_1$ is a large constant. Thus $\pi^*$ is
approximately dual to $\pi$ and is contained in $\G_N(C)$ for a suitable $C$. 
Two plates $\pi,\pi'$ have the same dual plate $\pi^*$ if and only if they
are parallel.

A {\em $\sigma$-cube} is a cube of side length $\sigma$ belonging to a suitable grid 
on $\R^{d+1}$; thus any two $\sigma$-cubes are either identical or have disjoint
interiors.  If $\sigma$ is fixed, for any $x\in\R^{d+1}$ we denote by $Q(x)$ the
$\sigma$-cube such that $x\in Q(x)$; should $x$ lie on boundary of two or more
cubes, we pick one of them arbitrarily.

If $R$ is a rectangular box (e.g. a tube or a plate), we will denote by
$cR$ the box obtained from $R$ by dilating it by a factor of $c$ 
about its center.

We define $\phi(x)=(1+|x|^2)^{-\frac{M}{2}}$ with $M$ large enough,
and $\phi_R=\phi\circ a_R^{-1}$, where $a_R$ is an affine map taking
the unit cube centered at $0$ to the rectangle $R$; thus $\phi_R$ is
roughly an indicator function of $R$ with ``Schwartz tails".  If ${\cal R}$
is a family of rectangular boxes (usually tubes or plates), we write
$\Phi_{\cal R}=\sum_{R\in{\cal R}}\phi_R$.
We note for future reference that if ${\cal R}$ is a family of separated
$\d$-plates or $\d$-tubes, then for any $\d$-cube $\Delta$ we have
\begin{equation}
\max_{\Delta}\Phi_{\cal R}\leq C\min_{\Delta}\Phi_{\cal R}
\label{harn}\end{equation}
where $C$ depends only on the choice of $M$.

We let $\psi(x):\R^{d+1}\rightarrow\R$ be a function such that 

\smallskip

1. $\psi=\eta^2$, where $\hat{\eta}$ is supported in a small ball centered at $0$.

2. $\psi\neq 0$ on a large cube centered at the origin.

3. The $\Z^{d+1}$ translations of $\psi$ form a partition of unity.

\smallskip\noindent
We also write $\psi_R=\psi\circ a_R^{-1}$ with $a_R$ as above. 

If a family of functions ${\cal F}_N$ is given for each $N$, we
will say that the functions in ${\cal F}_N$ are {\em essentially orthogonal} if
$$
\|\sum_{f\in {\cal F}_N} f\|_2^2\approx \sum_{f\in {\cal F}_N}\|f\|_2^2.
$$
For instance, functions
with finitely overlapping supports or Fourier supports are essentially
orthogonal.  Another important case will be discussed in the next section.

We conclude this section by stating and proving some basic properties of
the norm $\|\cdot\|_{p,mic}$.  

We first remark that, although we defined the norms $\|\ \|_{p,mic}$
and stated Theorem \ref{thm1} for functions with 
$\supp \widehat{f}\subset\G_N(1)$, we could as well have done so
for functions with $\supp \widehat{f}\subset\G_N(C)$ for some
constant $C$; furthermore, all of our estimates will continue
to hold in this case (modulo additional constant factors which will
be ignored in the sequel).  This is because $\G_N(C)$ can be covered 
by a bounded number of translates of sets of the form $\G_{N'}(1)$
with $N'\approx N$.

Observe also that for $p=\infty$ \bref{qa3} is just the trivial estimate
\begin{equation}\label{inf.trivial}
\|f\|_\infty\lesssim N^{\frac{d-1}{2}}\|f\|_{\infty,mic},
\end{equation}
which follows directly from the definition of $\|f\|_{\infty,mic}$ 
using that $f=\sum_\Theta \Xi_\Theta *f$ and there are at most
$N^{\frac{d-1}{2}}$ separated caps $\Theta$.

\begin{lemma}\label{interpol}
For all $p\geq 2$ we have
\begin{equation}
\|f\|_{p,mic}\lesssim\|f\|_2^{\frac{2}{p}}\|f\|_{\infty, mic}^{1-\frac{2}{p}}.
\label{looc5}\end{equation}
\end{lemma}

\underline{Proof} Let $f_\Theta=\Xi_\Theta*f$.  Since the supports
of $\widehat{f_\Theta} =y_\Theta\widehat{f}$ have finite overlap,
the functions $f_\Theta$ are essentially orthogonal:
\begin{equation}\label{www.98}
\|f\|_{2,mic}^2=\sum_{\Theta}\|f_\Theta\|^2_2
\approx \|f\|_2^2.
\end{equation}
Plugging this into the trivial estimate
\begin{eqnarray*}
\|f\|_{p,mic}^p &=& \sum_{\Theta}\|f_\Theta\|_p^p
\leq \max_{\Theta}\|f_\Theta\|_\infty^{p-2}
\cdot \sum_{\Theta}\|f_\Theta\|_2^2
\\
&=& \|f\|_{p,mic}^{p-2}\|f\|_{2,mic}^2
\end{eqnarray*}
we obtain the lemma.
\hfill$\square$

\bigskip

The next lemma describes the behaviour of $\|\cdot\|_{p,mic}$ under scaling.
If $f$ is a function on $\R^{d+1}$ and $R$ is a rectangle, we define
\[
T_Rf=(\psi_R f)\circ a_R=\psi\cdot(f\circ a_R).
\]

\begin{lemma}\label{Lemma 4.2}
Assume that $\supp \widehat{f}\subset\G_N(C)$, and let $Q$ be a $\sigma$-cube
for some $\d\lesssim\sigma\lesssim 1$. Then $\widehat{T_Qf}$ is supported
in $\G_{\sigma N'}(C')$, where $C'$ depends on $C$ and $N'\approx N$, and 
\begin{equation}
\|T_Qf\|_{p, mic}\lesssim \sigma^{-(\frac{d-1}{2}+\frac{1}{p})}
\|f\|_2^{\frac{2}{p}}\|f\|_{\infty, mic}^{1-\frac{2}{p}}
\label{ls20}\end{equation}
provided $p\geq 2$. Moreover, if $\sigma\geq N^{-\frac{1}{2}}$ then
\begin{equation}
\|T_Qf\|_2\lesssim N^{\frac{d-1}{4}}\|f\|_{\infty, mic}.
\label{ls30}\end{equation}
\end{lemma}

\underline{Proof}
We have $\widehat{T_Qf}=\widehat{\psi}*\widehat{(f \circ a_Q)}$.
By scaling, 
\begin{equation}\label{www.20}
\supp \widehat{f \circ a_Q}\subset \G_{\sigma N'}(\sigma C''),
\end{equation}
and the Fourier support statement follows since $\widehat{\psi}$
has compact support. 

Next, we prove \bref{ls20}. By \bref{looc5}, it suffices to do so for
$p=2$ and $p=\infty$.  For $p=2$, a standard argument using Schur's
test and \bref{www.20} shows that
\begin{equation}\label{Schur}
\|T_Q f\|_2\lesssim \sigma^{\frac{1}{2}}\|f\circ a_Q\|_2
=\sigma^{-\frac{d}{2}}\|f\|_2
\end{equation}
as required.  The $L^\infty$ bound follows from the fact that
a sector of angular length $(\sigma N)^{-\frac{1}{2}}$ intersects
$\lesssim \sigma^{-\frac{d-1}{2}}$ sectors of angular length
$N^{-\frac{1}{2}}$ (cf. \bref{inf.trivial}).

To prove \bref{ls30}, we write
\[
T_Qf=\sum_\Theta \psi \cdot((\Xi_{\Theta}\ast f)\circ a_Q)
\]
and observe that if $\sigma\geq N^{-\frac{1}{2}}$, 
the functions on the right are essentially orthogonal
since their Fourier supports have finite overlap. Hence
\begin{eqnarray*}
\|T_Q f\|_2^2&\lesssim&\sum_{\Theta}\|\psi\cdot((\Xi_{\Theta}\ast f)\circ a_Q)\|_2^2\\
&\lesssim&\sum_{\Theta}\|f\|_{\infty, mic}^2\cdot\|\psi\|_2^2\\
&\lesssim&N^{\frac{d-1}{2}}\|f\|_{\infty, mic}^2.
\end{eqnarray*}\hfill$\square$


\section{N-functions}\label{sec4}


\begin{definition}\label{Nfn}
An {\em $N$-function} is a function $f$ which has a decomposition
\[f=\sum_{\pi\in\P}f_{\pi},\]
where $\P=\P(f)$ is a separated family of $\d$-plates and
\begin{equation}\label{Nfna}
|f_{\pi}|\lesssim \phi_{\pi},
\end{equation}
\begin{equation}\label{Nfnb}
\supp \widehat{f_{\pi}}\subset \pi^*.
\end{equation}

\end{definition}

Such a decomposition is of course not unique; however, given an $N$-function
we will always fix a family of plates for $f$ and the associated functions
$f_{\pi}$. A {\em subfunction} of $f$ is a function of the form
\[
f_{\tilde{\P}}\stackrel{def}{=} \sum_{\pi\in\tilde{\P}}f_{\pi},
\]
where $\tilde{\P}$ is a subset of $\P$. 

$N$-functions are clearly Fourier supported in $\G_N(C)$.  Conversely,
Lemma \ref{Lemma 4.4} allows us to decompose functions $f$ with
$\supp \widehat{f}\subset \G_N(C)$ into $N$-functions while maintaining
control of the $\|\ \|_{p,mic}$ norms.

\begin{lemma}\label{Lemma 4.1}
Let $f$ be an $N$-function. Then we have the estimates
\begin{equation}\label{Nfn-linf}
\|f\|_\infty\lesssim N^{\frac{d-1}{2}},
\end{equation}
\begin{equation}
\|f\|_{p, mic}\lesssim \left(N^{-\frac{d+1}{2}}|\P|\right)^{\frac{1}{p}}
\hbox{ for }p\geq 2.
\label{ls3}\end{equation}
\end{lemma}

\underline{Proof} 
The estimate \bref{Nfn-linf} follows by counting the number
of $\d$-separated plates that can go through a fixed point.
It remains to prove \bref{ls3}. By \bref{looc5}, it suffices
to do so for $p=2$ and $p=\infty$. 

To check the case $p=2$ it suffices to verify that the functions
$f_{\pi}$ are essentially orthogonal.  Namely, we first write
$f=\sum_{\pi^*} g_{\pi^*}$, where $g_{\pi^*}$ is the sum of
those $f_\pi$ with $\pi$ dual to $\pi^*$; thus $g_{\pi^*}$ is
Fourier supported in $\pi^*$.  Since $\pi^*$ have bounded overlap,
$g_{\pi^*}$ are essentially orthogonal.  We may therefore assume
that all $\pi$ have the same dual plate $\pi^*$, hence are parallel
and have finite overlap.  But in this case it is easy to prove
essential orthogonality using the decay of $\phi_\pi$; the details
are left to the reader.

For $p=\infty$, we need to prove that
\[
\|f* \Xi_\Theta \|_\infty\lesssim 1
\]
for each $\Theta$.  We have 
$f* \Xi_{\Theta}=\sum_{\pi} f_\pi *\Xi_{\Theta}$,
where the only non-zero terms are those corresponding to
$\pi$ with $\pi^*\cap \G_{N,\Theta}(C)\neq\emptyset$.  But all such
$\pi$ are roughly parallel and therefore have finite overlap.
It follows that
\[
\|f* \Xi_\Theta \|_\infty
\lesssim \max_\pi \|f_\pi* \Xi_\Theta \|_\infty
\leq \int \phi_\pi \lesssim 1.
\]
$\ $\hfill$\square$

\begin{lemma}\label{Lemma 4.4}
Let $f$ be a function such that $\hat{f}$ is supported on $\Gamma_N(C)$
and $\|f\|_{\infty,mic}<\infty$. Then there are $N$-functions
$f_{\l}$, with dyadic $\l$ satisfying
\begin{equation}
\l\lesssim \|f\|_{\infty, mic},
\label{katr1}\end{equation}
such that 
\begin{equation}\label{Nfn-dec}
f\approx \sum_{\l}\l f_{\l},
\end{equation}
\begin{equation}
\sum_{\l}\l^p\d^{\frac{d+1}{2}}|\P(f_{\l})|\lesssim\|f\|_{p, mic}^p
\label{wa1}\end{equation}
for each fixed $p\in [2,\infty)$.
\end{lemma}

\underline{Proof} 
We may assume that $\hat f$ is supported in $\G_{N,\Theta}(C)$
for some $\Theta$, so that $\|f\|_{p,mic}=\|f\|_p$.  
We fix a plate $\pi$ so that $\supp \hat f\subset \frac{1}{2}\pi^*$.
This is possible if the constant in the definition of $\pi^*$ was
chosen large enough.  
Let $\{\pi_j\}$ be a tiling of $\R^{d+1}$ by translates of $\pi$,
and let $\psi_j=\psi_{\pi_j}$.  Define
\[
\P_{\l}=\{\pi_j:\ \l\leq \|\psi_{j}f\|_\infty\leq 2\l\},
\]
and
\[
f_\l=\sum_{\pi_j\in\P_\l}\l^{-1}\psi_{j}^2 f.
\]
Clearly, $\P_{\l}$ can be non-empty only for $\l$ as in
\bref{katr1}.  

We first show that each $f_{\l}$ is an $N$-function.  Fix 
a plate $\pi_j\in\P_\l$ and let $f_{j}:=\l^{-1}\psi_{j}^2 f$.
It is clear from the definition that $f_j$ satisfies \bref{Nfna}.
\bref{Nfnb} follows from the fact that $\widehat{f_j}=
\widehat{\psi_j^2}*\hat f$ and $\widehat{\psi_j^2}$ is
supported in $c\pi^*_0$, where $c\ll 1$ is a small constant
and $\pi^*_0$ is a translate of $\pi^*$ centered at $0$.

Next, we have
\begin{equation}\label{wa-partition}
1\lesssim \sum_j \psi^2_j \lesssim \sum_j \psi_j=1.
\end{equation}
Since $\sum_\l \l f_\l=\sum_j \psi_j^2 f$, \bref{Nfn-dec} follows.
It remains to prove \bref{wa1}. By Bernstein's inequality
and \bref{wa-partition},
\[
\l^p\approx \|\psi_jf\|^p_{\infty}\lesssim |\pi^*|\|\psi_jf\|_p^p
\approx \d^{-\frac{d+1}{2}}\|\psi_jf\|_p^p.
\]
Hence
\[
\sum_{\l}\l^p\d^{\frac{d+1}{2}}|\P_{\l}|
\lesssim \sum_j\|\psi_jf\|_p^p\lesssim\|f\|_p^p
\]
as claimed.\hfill$\square$


\section{Proof of Theorem \ref{thm1}}\label{sec5}


The main step in the proof of Theorem \ref{thm1} is the following 
inductive argument.

\begin{definition}\label{def-P}
We say that $P(p,\alpha)$ holds if for all functions $f$ such that
$\supp\widehat{f}\subset \G_N(1)$ and $\|f\|_{\infty, mic}\leq 1$ we have
\begin{equation}
|\{ |f|>\l\}|\lesssim \l^{-p} \d^{d-\frac{(d-1)p}{2}-\alpha}\|f\|_2^2\ ,
\label{4.1}\end{equation}
provided that $\d$ is small enough. 
\end{definition}

\begin{proposition}\label{scales}
Fix $p>p_d$ and suppose that $P(p,\alpha)$ holds. 
Then $P(p,\beta)$ holds for any $\beta>(1-\frac{\e_0}{4})\alpha$. 
\end{proposition}

In this section we will prove Theorem \ref{thm1} assuming that Proposition
\ref{scales} is known.  The remaining sections will be devoted to the proof
of Proposition \ref{scales}.

\begin{corollary}\label{cor-iterate}
Assume that we have proved Proposition \ref{scales}.  Then $P(p,\alpha)$
holds for all $\alpha>0$.
\end{corollary}

\underline{Proof} Let $f$ satisfy the assumptions of Definition 
\ref{def-P}.   It suffices to show that \bref{4.1} holds for some
(large) $\alpha>0$, since then the conclusion will follow by iterating
Proposition \ref{scales}.

The left side of \bref{4.1} can be non-zero only when
\[
\l\leq \|f\|_\infty\lesssim \d^{-\frac{d-1}{2}},
\]
where the last inequality follows from \bref{inf.trivial}.  On the other
hand, \bref{4.1} follows from Tchebyshev's inequality if
\begin{equation}\label{cheb}
\l^{p-2}\lesssim  \d^{d-\frac{(d-1)p}{2}-\alpha}.
\end{equation}
But this holds for all $\l$ as above if $\alpha$ has been chosen large enough.
\qed

\begin{lemma}\label{strong-type}
If \bref{4.1} holds, then the corresponding strong type estimate
\begin{equation}
\|f\|_p^p\lessapprox  \d^{d-\frac{(d-1)p}{2}-\alpha} \|f\|_2^2
\label{mnmnmn}\end{equation}
also holds for the same class of $f$.
\end{lemma}

\underline{Proof}
Write $|f|\approx \sum_\l \l\chi_{|f|\approx\l}$ with dyadic $\l$.
By \bref{inf.trivial} we have $\|f\|_\infty\lesssim \d^{-\frac{d-1}{2}}$;
on the other hand, \bref{mnmnmn} is trivial for functions with
$\|f\|_\infty\lesssim\d^{n}$ for $n$ large enough.  Therefore
the lemma follows by summing \bref{4.1} over dyadic $\l$ 
with $\d^n\lesssim \l\lesssim \d^{-\frac{d-1}{2}}$.
\qed

~

\underline{Proof of Theorem \ref{thm1}, given Proposition \ref{scales}.}
Let $p>p_d$.
Let $f$ be Fourier supported in $\Gamma_N(1)$, and assume that
\begin{equation}\label{www.99}
\|f\|_{p, mic}\leq 1.
\end{equation}
Fix $\e>0$. We will prove that
\begin{equation}\label{www.100}
\|f\|_{L^p(Q_0)} \lesssim\d^{\frac{d}{p}-\frac{d-1}{2}-\e},
\end{equation}
where $Q_0$ is the unit cube.  A standard argument using the partition
of unity $\{\psi(x-j)\}_{j\in\Z^{d+1}}$ will then yield Theorem
\bref{thm1}; the details are left to the reader.

Observe also that \bref{www.99} implies that $\|f\|_\infty\lesssim \d^{-K}$
for some large constant $K$.  Indeed, let $f_\Theta=\Xi_\Theta *f$.
Since $\widehat{f_\Theta}$ is supported in $\G_{N,\Theta}(C)$, we may
apply Bernstein's inequality to $f_\Theta$, obtaining
\[
\|f_\Theta\|_\infty\lesssim |\G_{N,\Theta}(C)|^\frac{1}{p}\|f_\Theta\|_p
\lesssim N^{\frac{d+1}{2p}}.
\]
Summing over $\Theta$, we obtain the claimed bound.

Let $f\approx \sum_{\l\leq \d^{-K}}\l f_{\l}$ be the decomposition
given by Lemma \ref{Lemma 4.4}.  By \bref{wa1} and Lemma \ref{Lemma 4.1},
the $N$-functions $f_\l$ satisfy
\[
\|f_{\l}\|_{\infty, mic}\lesssim 1
\]
and (using also \bref{www.98})
\begin{equation}\label{www.97}
\|f_{\l}\|_2^2\approx \|f_{\l}\|_{2,mic}^2\lesssim\l^{-p}.
\end{equation}

By Corollary \ref{cor-iterate} and Lemma \ref{strong-type}, we have
\bref{mnmnmn} for any $\alpha>0$.  Substituting \bref{www.97} in
\bref{mnmnmn}, we obtain
\[
\|f_{\l}\|^p_p\lesssim\l^{-p} \d^{d-\frac{(d-1)p}{2}-p\e},
\]
which is \bref{www.100}.

Note also that $\sum_{\l\leq\d^{K'}}\l f_{\l}$ with $K'$ large enough is bounded
pointwise by a large power of $\d$, hence satisfies \bref{www.100}.
Since there are logarithmically many dyadic $\l$ with 
$\d^{K'}\lesssim\l\lesssim\d^{-K}$, the result follows by summing over $\l$.
\qed


\section{A few geometrical lemmas}\label{sec2}


In this section we collect the geometrical information needed in the proof
of Proposition \ref{scales}.

We begin with some preliminaries.  Let $\P$ be a family of $\d$-plates. 
For each $\pi\in\P$ we pick a $\d^{\frac{1}{2}}$-tube $\tau$ containing
$\pi$, thus obtaining a family of tubes $\tilde\T$.  We let $\T$ be
a maximal $\sqrt{\d}$-separated subset of $\tilde \T$.  Replacing the 
tubes $\tau\in\T$ by their dilates $C\tau$ if necessary, we obtain
a separated family of $C\sqrt{\d}$-tubes $\T(\P)$ so that each $\pi\in\P$
is contained in some $\tau(\pi)\in\T(\P)$ (if $\pi$ is contained
in more than one tube in $\T(\P)$, we pick one arbitrarily),
and that the plates $\pi$ with the same $\tau(\pi)=\tau$ are all parallel.  

For each $\tau\in\T(\P)$ we define
\[
X_{\P}(\tau)=\{\pi\in\P:\tau(\pi)=\tau\}.
\]
Thus $\P$ is the disjoint union of the $X_{\P}(\tau)$'s where $\tau$ runs
over $\T(\P)$. 
The {\em type $r$ component} of $\P$ is the subset $\P_r\subset\P$ defined by
\[
\P_r=\bigcup_{ r\leq |X_{\P}(\tau)|\leq 2r} X_{\P}(\tau).
\]
We will say that $\P$ is {\em type $r$} if $\P=\P_r$.

\begin{lemma}\label{Lemma 2.9}
Assume that $\P$ is type $r$. Then for any
$\d^{\frac{1}{2}}$-cube $Q$ we have
\[
\int\Phi_{\P}\phi_Q\lesssim\d^{\frac{1}{2}}r\int\Phi_{\T(\P)}\phi_Q.
\]
\end{lemma}

\underline{Proof} We first prove that for any $\d^{\frac{1}{2}}$-cube $Q$
and for each $\tau\in\T(\P)$
\[
\int_Q \sum_{\pi:\tau(\pi)=\tau} \phi_{\pi}\lesssim\d^{\frac{1}{2}}r\int_Q\phi_{\tau}.
\]
Let $T$ be the infinite tube extending $\tau$ in the direction of its
longest axis.  If $Q\cap 100T\neq\emptyset$, we have
\[
\int_Q \sum_{\pi:\tau(\pi)=\tau} \phi_{\pi}
\lesssim \d^{\half}r|Q| \lesssim\d^{\frac{1}{2}}r\int_Q\phi_{\tau}.
\]
If on the other hand $Q\cap 100T =\emptyset$, we have the pointwise
estimate 
$\sum_{\pi:\tau(\pi)=\tau} \phi_{\pi}(x)\lesssim\d^{\frac{1}{2}}r\phi_{\tau}(x)$.
This proves our claim.

It follows that
\[
\int_Q\Phi_{\P}\lesssim\d^{\frac{1}{2}}r\int_Q\Phi_{\T(\P)}
\]
for each $\d^{\half}$-cube $Q$.  The lemma is now easily proved by writing
$\phi_Q=\sum_{Q'}\phi_Q\chi_{Q'}$ and using estimates like \bref{harn}.
\qed

~

The main geometrical argument is as follows.

\begin{lemma}\label{Lemma 2.3}
Let $\P$ be a separated family of $\d$-plates, and let $W\subset
\R^{d+1}$.  Assume that $t\approx\d^{c}$ with $0<c<1$.  Then there
is a relation $\sim$ between plates in $P$ and $t$-cubes $Q$ such that
\begin{equation}\label{propR}
|\{Q:\ \pi\sim Q\}|\lessapprox 1\hbox{ for all }\pi\in\P
\end{equation}
(with the implicit constants independent of $\pi$), and
\begin{equation}\label{badI}
I_b\lessapprox t^{-3d}|W||\P|^{\frac{1}{2}},
\end{equation}
where
\[
I_b=\int_W \sum_{\pi\in\P,\pi\not\sim Q(x)}\chi_\pi(x)
=\sum_{\pi\in\P}|\{x\in W\cap\pi:\ Q(x)\not\sim \pi\}.
\]
\end{lemma}

\underline{Proof} We may assume that all the plates and points are contained
in a large cube of side length $\lesssim 1$.

We define a relation $\sim$ by declaring that each $\pi\in\P$ is related to
$Q(\pi)$ and its neighbours, where $Q(\pi)$ is the $t$-cube with maximal
$|W\cap Q\cap\pi|$ (if there is more than one such cube, choose one
arbitrarily).  It is clear that \bref{propR} holds; we need to prove \bref{badI}.

We can pigeonhole to obtain $\P'\subset \P$ and $\nu$ so that
\begin{equation}\label{inc.yikes}
|I_b|\lessapprox \nu |\P'|
\end{equation}
and 
\[
|\{x\in W\cap\pi:\ \pi\not\sim Q(x)\}|\approx\nu\hbox{ for each }\pi\in \P'.
\]
Hence for each $\pi\in\P'$ there is a cube $Q'(\pi)$ such that $\pi\not\sim
Q'(\pi)$ and $|W\cap Q'(\pi) \cap \pi|\gtrsim t\nu$.
By the choice of $Q(\pi)$, we also have $|W\cap Q(\pi)\cap \pi|\gtrsim t\nu$.
Since the number of all possible pairs of $t$-cubes $(Q,Q')$ is
$\lesssim t^{-2d-2}$, there are cubes $Q,Q'$ such that $Q=Q(\pi)$ and
$Q'=Q'(\pi)$ for at least $t^{2d+2}|\P'|$ plates $\pi\in\P'$.
Fix such $Q$ and $Q'$, and consider the quantity
\[
A=\sum_{\pi\in\P'}|W\cap Q\cap \pi|\cdot |W\cap Q'\cap\pi|.
\]
From the above estimates we have
\[
A\gtrsim t^{2d+2}|\P'|\cdot (t\nu)^2=t^{2d+4}\nu^2|\P'|.
\]
On the other hand, we can rewrite $A$ as
\[
A=\int_{W\cap Q}\int_{W\cap Q'}\sum_{\pi\in\P'}\chi_\pi(x)\chi_\pi(x')dx'dx.
\]
Since $Q'$ is at distance at least $t$ from $Q$,
for any $x\in Q,x'\in Q'$ there are $\lesssim t^{-d+1}$ separated 
plates that go through both $x$ and $x'$; in other words, 
$\sum_{\pi\in\P'}\chi_\pi(x)\chi_\pi(x')\lesssim t^{-d+1}$.  Hence
\[
A\lesssim t^{-d+1}|W\cap Q|\cdot|W\cap Q'|
\lesssim t^{-d+1}|W|^2.
\]
Comparing the upper and lower bounds for $A$, we find that
\[
\nu\lesssim t^{-\frac{3d+3}{2}}|W|\,|\P'|^{-1/2}.
\]
Plugging this into \bref{inc.yikes} and using that $|\P'|\leq|\P|$,
we obtain \bref{badI}.
\hfill$\square$

~


We now prove a version of Lemma \ref{Lemma 2.3} with ``Schwartz tails";
the argument here is quite standard.

\begin{lemma}\label{Lemma 2.7} Let $W\subset\R^{d+1}$ and let $\P$ be
a family of separated $\d$-plates. Fix a large constant $M_0$. Then, if
the constant $M$ in the definition of $\phi$ has been chosen large
enough, there is a relation between $t$-cubes and 
plates in $\P$ satisfying \bref{propR} and such that
\begin{equation}
\int_W\Phi_{\P}^b\lessapprox t^{-3d}|\P|^{\frac{1}{2}}|W| +\d^{M_0}|W|,
\label{vu2}\end{equation}
where
\begin{equation}
\Phi^b_{\P}(x)=\sum_{\pi\in\P,\pi\not\sim Q(x)}\phi_{\pi}(x)
\label{sch2}\end{equation} 
\end{lemma}

\underline{Proof} 
For each $\pi\in\P$ we have a decomposition $\R^{d+1}=\bigcup_{k=0}^\infty
2^k\pi$, with $\phi_\pi\approx 2^{-kM}$ on $2^{k+1}\pi\setminus 2^k\pi$.
Therefore 
\begin{equation}\label{www.1}
\int_W \Phi^b_{\P} 
\lesssim\int_W \sum_{\pi\in\P,\pi\not\sim Q(x)}\sum_{k=0}^\infty 2^{-kM}\chi_{2^k\pi}(x)dx.
\end{equation}
Pick $K$ so that $c_0\log(t\d^{-\half})\leq K\leq 2c_0 \log(t\d^{-\half})$
for a sufficiently small constant $c_0>0$. Then the last integral is bounded by
\begin{equation}\label{www.2}
\sum_{k=0}^K 2^{-kM} \int_W \sum_{\pi\in\P,\pi\not\sim Q(x)}\chi_{2^k\pi}(x)dx
+\int_W \sum_{k>K} \sum_{\pi\in\P} 2^{-kM}\chi_{2^k\pi}(x)dx.
\end{equation}

Fix a value of $k\leq K$.  Let $\P_k$ be a maximal subset of $\P$ such that the
plates $\{2^k\pi:\ \pi\in\P_k\}$ are separated (on scale $2^k\d$).  Then for
each $\pi\in\P$ there is a $\pi_k\in\P_k$ such that $2^k\pi$ is comparable
to $2^k\pi_k$. Let $c$ be a constant such that $2^k\pi\subset 2^kc\pi_k$
for all $\pi\in\P$.

If $c_0$ was chosen suitably small, $2^K \sqrt{\d}$ is sufficiently small
compared to $t$, so that we may apply Lemma \ref{Lemma 2.3} to the plate
family $\{2^kc\pi: \pi\in\P_k\}$. The relation thus obtained will
be denoted by $\sim_{k,0}$.
We now define a relation $\sim_k$ between plates in $\P$ and $t$-cubes
by declaring that 
\[
\pi\sim_k Q\ \Leftrightarrow\ 2^k c\pi_k\sim_{k,0}Q.
\]
Then for each $\pi$ we have 
\[
\{x\in W\cap 2^k\pi:\ \pi\not\sim_k Q(x)\}\subset
\{x\in W\cap 2^kc\pi_k:\ 2^k c\pi_k\not\sim_{k,0} Q(x)\}.
\]
Since there are at most $\lesssim 2^{k(d-1)}$ plates $\pi$ with the same
$\pi_k$, we have from \bref{badI}
\begin{eqnarray*}
\int_W \sum_{\pi\in\P,\pi\not\sim Q(x)}\chi_{2^k\pi}(x)dx
& = &\sum_{\pi\in\P}|\{x\in W\cap 2^k\pi:\ \pi\not\sim_k Q(x)\}|
\\
& \lesssim & 2^{k(d-1)}\sum_{\pi\in\P_k}
|\{x\in W\cap 2^kc\pi:\ 2^kc\pi\not\sim_{k,0} Q(x)\}|
\\
& \lessapprox & 2^{k(d-1)}t^{-3d}|W|\,|\P_k|^\half
\\
&\lessapprox & 2^{k(d-1)}t^{-3d}|W|\,|\P|^\half.
\end{eqnarray*}

Finally, we define a relation $\sim$ as follows:
\[
\pi\sim Q\ \Leftrightarrow\ \pi\sim_k Q\hbox{ for some }k\leq K.
\]
We have $K\lessapprox  1$, hence \bref{propR} holds.  Moreover, 
the first term in \bref{www.2} is bounded by
\begin{eqnarray*}
&\lesssim & \sum_{k=0}^K 2^{-kM} \int_W \sum_{\pi\not\sim_k Q(x)}\chi_{2^k\pi}(x)dx
\\
&\lessapprox& \sum_{k=0}^K 2^{-kM}2^{k(d-1)}t^{-3d}|W|\,|\P|^\half
\\
&\lessapprox& t^{-3d}|W|\,|\P|^\half.
\end{eqnarray*}
It remains to estimate the second term in \bref{www.2}. For any $x$ we have
\[
|\{\pi\in \P:\ x\in 2^k\pi\}|\lesssim 2^{k(d+1)}\d^{-2d}.
\]
Hence, if $M$ was chosen large enough, we obtain that
\begin{eqnarray*}
\int_W \sum_{k>K} \sum_{\pi\in\P} 2^{-kM}\chi_{2^k\pi}(x)dx
&\lesssim& \int_W \sum_{k>K} 2^{-kM+k(d+1)}\d^{-2d}dx\\
&\lesssim& 2^{-KM/2}\d^{-2d}|W|\\
&\lesssim& \d^{M_0}|W|,
\end{eqnarray*}
where at the last step we used that $K\gtrsim \log(\frac{1}{\d})$.
\hfill$\square$

~

We will also need the analogue of Lemma \ref{Lemma 2.7} with plates replaced by tubes.
The proof is identical to that of Lemma \ref{Lemma 2.7}, therefore we omit it.

\begin{lemma}\label{Lemma 2.8}
Let $W\subset\R^{d+1}$ and let $\T$ be a family of separated $\sqrt{\d}$-tubes.
Fix a large constant $M_0$. Then there is a relation between $t$-cubes and 
tubes in $\T$ satisfying \bref{propR} and such that
\begin{equation}
\int_W\Phi_{\T}^b\lessapprox t^{-3d}|\T|^{\frac{1}{2}}|W| +\d^{M_0}|W|,
\label{vu2a}\end{equation}
where
\begin{equation}
\Phi^b_{\T}(x)=\sum_{\tau\in\T,\tau\not\sim Q(x)}\phi_{\tau}(x)
\label{sch2a}\end{equation} 
\end{lemma}


\section{A localization property}\label{sec3}


In this section and throughout the rest of this paper we will always assume
that $t$ is a dyadic nuumber such that $t\approx \d^{\e_0}$ with $\e_0>0$
sufficiently small.

\begin{definition}\label{localize}
Let $f$ be an $N$-function. We say that $f$ {\em localizes} at $\lambda$
if there are subfunctions $f^Q$ of $f$, where $Q$ runs over $t$-cubes, such that
\begin{equation}
\sum_Q|\P(f^Q)|\lessapprox |\P(f)|
\label{vo9}\end{equation}
and
\begin{equation}\label{zzz.45} 
|\{|f|\geq\l\}|\lessapprox  \sum_Q |Q\cap \{|f^Q|\gtrapprox \l\}|.
\end{equation}
\end{definition}

\begin{lemma}\label{Lemma 3.1}
Let $f$ be an $N$-function with plate family $\P$. Assume that 
\begin{equation}
|\P|\leq t^{4d}\l^2.
\label{vo8}\end{equation} 
Then  $f$ localizes at $\l$.
\end{lemma}

\underline{Proof of Lemma \ref{Lemma 3.1}} Let $k=|\P|$; since $k\geq 1$,
\bref{vo8} implies that $\l\geq 1$. 
Let also $W=\{|f|\geq\l\}$.
Since $|f_{\pi}|\leq\phi_{\pi}$ and $\P$ is separated we  also have
$\l\leq\Phi_{\P}(x)\lesssim\d^{-\frac{d-1}{2}}$ if $x\in W$.

We apply Lemma \ref{Lemma 2.7} to $\P$ and $W$, obtaining a relation $\sim$ such that
\[
\int_{W}\Phi_{\P}^b\lessapprox t^{-3d }k^{\frac{1}{2}}|W|
\lesssim t\l|W|,
\]
where the last inequality follows from \bref{vo8}.  Hence there is
a subset $W^*\subset W$ with proportional measure such that
\begin{equation}
\Phi_{\P}^b(x)\lesssim t\l\hbox{ for }x\in W^*.
\label{vu8}\end{equation}
Define for each $Q$
\[
f^Q=\sum_{\pi\sim Q}f_\pi.
\]
Then for $x\in W^*\cap Q$ we have
\[
|f(x)-f^Q(x)|=|\sum_{\pi\not\sim Q}f_\pi(x)|
\lesssim \Phi_{\P}^b(x)\lesssim t\l,
\]
so that $|f^Q|\gtrsim\l$ on $W^*\cap Q$.  It remains only
to observe that the bound \bref{vo9} follows from 
\bref{propR}. \hfill$\square$

~

Theorem \bref{thm1} with $p_d=2+\frac{8}{d-3}$ can be proved using only
Lemma \ref{Lemma 3.1}.  However, this does not give any result for $d=3$.  
We therefore prove a similar lemma with the assumption \bref{vo8} replaced
by \bref{vo5}, thus gaining an additional factor of nearly $\d^{-\frac{d-1}{4}}$ 
when $\l$ is close to its maximum possible value $\d^{-\frac{d-1}{2}}$.
The conclusion of the lemma is somewhat weaker: essentially, it will
allow us to obtain a localization effect on one of the two scales
$N$ or $\sqrt{N}$.

\begin{lemma}\label{Lemma 3.2} 
Let $f$ be an $N$-function and assume that 
\begin{equation}|\P(f)|\leq t^{20d}\d^{\frac{3d-3}{4}}\l^4.
\label{vo5}\end{equation} 
Then either $f$ localizes at $\l$, or else there is a subfunction
$f^*$ of $f$ such that $|f^*|\gtrapprox \l$ on a logarithmic fraction
of $\{|f|\geq\l\}$ $W$, and 
\begin{equation}
\|\psi_{\Delta}f^*\|_{2}^2\lesssim t^{5d} \d^{\frac{3}{4}(d+1)}\l^2
\label{vo5.1}\end{equation}
for each $\d^{\frac{1}{2}}$-cube $\Delta$.

\end{lemma}

\underline{Proof} We note that $\l>1$, and let $k=|\P|$.
Let $\P_r$ be the type $r$ component of $\P$, then for
some $r$ we must have 
\[
|W|\gtrapprox |\{|f|\geq\l\}|,
\]
where $W=\{x:\ |f(x)|\geq\l,\ |f_{\P_r}(x)|\gtrapprox \l\}$.
With this value of $r$, we let $\T=\T(\P_r)$ be the family of
$\sqrt{\d}$-tubes defined in Section \ref{sec2}.
We clearly have
\[
|\T|\lessapprox kr^{-1}.
\]
We now consider two cases. 

~

{\it Case 1:} $\l\geq t^{-4d}(\frac{k}{r})^{\frac{1}{2}}$.
We claim that $f$ localizes; the proof is similar to that of
Lemma \ref{Lemma 3.1},  except that we use Lemma \ref{Lemma 2.8}
instead of Lemma \ref{Lemma 2.7}.
By Lemma \ref{Lemma 2.8}, there is a relation 
$\sim$ between tubes from $\T$ and $t$-cubes satisfying \bref{propR} and
\[
\int_{W}\Phi_{\T}^b
\lessapprox t^{-3d} \Big(\frac{k}{r}\Big)^{\frac{1}{2}}|W|
\lesssim t\l |W|.
\]
Hence there is a subset $W^*$ of $W$ with proportional measure such
that $\Phi_{\T}^b\lesssim t\l$ on $W^*$.
We define a relation between plates $\pi\in\P_r$ and $t$-cubes via
$\pi\sim Q$ if $\tau(\pi)\sim Q$, and let
\[
f^Q=\sum_{\pi\sim Q}f_{\pi}.
\]
Then for $x\in W^*\cap Q$ we have
\[
|f(x)-f^Q(x)|=|\sum_{\pi\not\sim Q}f_\pi(x)|
\lesssim \Phi_{\T}^b(x)\lesssim t\l,
\]
hence $|f^Q|\gtrsim\l$ on $W^*\cap Q$.  
The bound \bref{vo9} follows from \bref{propR}. \hfill$\square$

~ 

{\it Case 2:} $\l\leq t^{-4d}(\frac{k}{r})^{\frac{1}{2}}$.
We will show that in this case $f_{\P_r}$ satisfies \bref{vo5.1}.
Fix a $\d^{\half}$-cube $\Delta$. A slight modification of the
argument in the proof of Lemma \ref{Lemma 4.1} shows that the
functions $\psi_{\Delta}f_{\pi}$ are essentially orthogonal.
Hence
\[
\|\psi_{\Delta}f_{\P_r}\|_2^2\approx 
\|\sum_{\pi\in\P_r}\psi_{\Delta}f_{\pi}\|_2^2
\lesssim\int\sum_{\pi\in\P_r}|f_{\pi}|^2\phi_{\Delta}
\lesssim\int\Phi_{\P_r}\phi_{\Delta}.
\]
By Lemma \ref{Lemma 2.9}, the last integral is 
bounded by $\d^{\frac{1}{2}}r\int\Phi_{\T}\phi_{\Delta}$.  
An easy Schwartz tails calculation using the pointwise bound 
$\Phi_{\T}\lesssim \d^{-\frac{d-1}{2}}$ shows that
\[
\int\Phi_{\T}\phi_{\Delta}\lesssim\d^{-\frac{d-1}{2}}\d^{\frac{d+1}{2}}=\d.
\]
We are also assuming that $r\leq t^{-8d} k\l^{-2}$.  Therefore
\[
\|\psi_{\Delta}f_{\P_r}\|_2^2
\lesssim \d^{\frac{3}{2}}r
\lesssim t^{-8d}\d^{\frac{3}{2}}\l^{-2}k
\lesssim t^{5d}\d^{\frac{3d+3}{4}}\l^2,
\]
where at the last step we used \bref{vo5}.
\hfill$\square$


\section{Proof of Proposition \ref{scales}}\label{sec5a}


In this section we will prove Proposition \ref{scales}.  The general
scheme of the proof is as follows.  We will see in Lemma \ref{Lemma 5.2}
that it is easy to prove the proposition for $N$-functions which localize; 
therefore the main issue is to obtain the localization effect on some
suitable scale $\tilde N$.   

In Lemma \ref{Lemma 5.1} we decompose $f$ into functions with Fourier
support in cone sectors of size roughly $N\times N^{\frac{3}{4}}\times
\dots\times N^{\frac{3}{4}}$, use a Lorentz transformation to rescale
each of these sectors to a neighbourhood of $\G_{\sqrt{N}}$, and
apply the inductive hypothesis to the functions thus obtained.  
We are then left with the task of estimating the measure of the sets
of large values of certain parts of $f$ which can be rescaled to
$\sqrt{N}$-functions $f_\Delta$, with good estimates on the cardinalities of
the corresponding plate families.
Namely, if $p>p_d$ then it can be shown that $f_\Delta$ satisfy the
assumptions of either Lemma \ref{Lemma 3.1} or Lemma \ref{Lemma 3.2}.
Thus either $f_\Delta$ localize, or else we can change scales again
and obtain localization for an appropriate further decomposition of
$f_\Delta$.  Applying the inductive hypothesis again on scale
slightly smaller than $\sqrt{N}$, we obtain the proposition.

\begin{lemma}\label{Lemma 5.1}  Assume that $P(p,\alpha)$ is known
for some $p$ and $\alpha$.
Let $f$ be Fourier-supported in $\Gamma_N(C)$  and such
that $\|f\|_{\infty, mic}\lesssim 1$. Then for any $\l\geq\d^{C_9}$ 
and for any $\e>0$ there is a $\l_*$ and a collection 
of $\sqrt{N}$-functions $\{f_{\Delta}\}$  so that
a logarithmic fraction of $\{|f|\geq\l\}$ is contained in
$\bigcup_\Delta a_\Delta^{-1}(\{|f_\Delta|\geq\l_*\})$, and
\begin{equation}\label{vp.lambda}
\l \d^{\frac{d-1}{4}+\e}\lesssim \l_*\lesssim \d^{-\frac{d-1}{4}},
\end{equation}
\begin{equation}
|\P(f_{\Delta})|\leq\d^{-C\e}(\frac{\l_*}{\l})^2\d^{-\frac{3d+3}{4}}
\|\psi_\Delta f\|_2^2,
\label{vp2}\end{equation}
\begin{equation}
\sum_{\Delta}|\P(f_{\Delta})|\leq\d^{-C\e}
(\frac{\l_*}{\l})^{p}\d^{-\frac{3d+3}{4}}
\d^{\frac{d}{2}-\frac{(d-1)p}{4}-\frac{\alpha}{2}} \|f\|_2^2.
\label{vp3}\end{equation}
\end{lemma}

\underline{Proof} We write $f=\sum_\Delta \psi_\Delta f$, where
$\Delta $ runs over $\sqrt{\delta}$-cubes.  Fix a small $\e>0$.
It is an easy exercise to prove that
\begin{equation}
\{|f|\geq\l\}\subset\bigcup_\Delta \{|\psi_\Delta f|\geq c \d^{\e}\l\},
\label{ls13}\end{equation}
using \bref{inf.trivial} and the Schwartz decay of $\psi$.

For each $\Delta$ we apply Lemma \ref{Lemma 4.4} at scale $\sqrt{N}$
to $T_\Delta f$, obtaining a decomposition
\[
T_\Delta f \approx \sum_{h}hg^\Delta_h,
\]
where $g^\Delta_h$ are $\sqrt{N}$-functions.  By \bref{katr1} and \bref{ls20} we have
\begin{equation}\label{ls21}
h\lesssim\|T_\Delta f\|_{\infty, mic}
\lesssim \d^{-\frac{d-1}{4}}.
\end{equation}
Since we also assume that $\l\geq\d^{C_9}$, there are logarithmically many 
relevant dyadic values of $h$.  We may therefore choose $h=h(\Delta)$ so that 
a logarithmic fraction of $\{|T_\Delta f|\geq\d^{\e}\l\}$ 
is contained in the set $\{|hg^\Delta_h|\geq\d^{2\e}\l\}$. 
Finally, we pigeonhole to get a value of $h$ so that
a logarithmic fraction of $\{|f|\geq \l\}$ is contained in
$\bigcup_\Delta a_\Delta^{-1}(\{|hg^\Delta_h|\geq\d^{2\e}\l\})$.

Let $\l_*=\d^{2\e}\l h^{-1}$ and $f_\Delta=g^\Delta_h$,
with this value of $h$.  The lower bound in \bref{vp.lambda}
follows from \bref{ls21}.  To obtain the upper bound, we use that
\[
\l_*=\d^{2\e}\l h^{-1}\leq \|g_h^\Delta\|_{\infty}
\lesssim N^{\frac{d-1}{4}},
\]
where the last inequality follows from \bref{Nfn-linf}.

Let $\P_\Delta$ be the plate family for $f_\Delta$.
Applying \bref{wa1} on scale $\sqrt{N}$ and using that $h=\d^{2\e}\frac{\l}{\l_*}$ we obtain
\begin{equation}
|\P_\Delta|\lesssim N^{\frac{d+1}{4}}(\d^{2\e}\frac{\l}{\l_*})^{-p}
\|T_\Delta f\|_{p, mic}^p.
\label{ls14}\end{equation}
Applying \bref{ls14} with $p=2$ and then using that
$\|T_\Delta f\|^2_2=\d^{-(d+1)/2}\|\psi_\Delta f\|_2^2$,
we obtain \bref{vp2}.

It remains to prove \bref{vp3}.  By \bref{ls14}, it suffices to show that
\begin{equation}\label{zzz.10}
\sum_\Delta \|f_\Delta\|_{p,mic}^p\lesssim \d^{-\frac{d+1}{2}-C\e}
\d^{\frac{d}{2}-\frac{(d-1)p}{4}-\frac{\alpha}{2}}\|f\|_2^2.
\end{equation}

Let $\{\Psi\}$ be a finitely overlapping covering of $S^{d-1}$ by
spherical caps of angular
length $N^{-1/4}$, so that the corresponding cone sectors $\G_{N,\Psi}(C)$
have dimensions roughly $CN\times CN^{\frac{3}{4}}\times\dots\times
CN^{\frac{3}{4}}\times C$. 
Let also $\Xi_{\Psi}$ are functions whose Fourier transforms
agree on $\G_N(C)$ with a partition of unity subordinate to the
covering $\{\G_{N,\Psi}(C)\}$ of $\G_N(C)$.

We first claim that
\begin{equation}\label{ls10}
\sum_\Delta \|f_\Delta \|_{p, mic}^p
\lesssim \d^{-\frac{d+1}{2}}\sum_{\Psi}\|\Xi_{\Psi}\ast f\|_p^p,
\end{equation}
Indeed, using the definition of $\|\cdot\|_{p,mic}$ and
rescaling $x\to \d^{-\frac{d+1}{2}}x$ we obtain that
\[
\sum_\Delta \|f_\Delta \|_{p, mic}^p =
\d^{\frac{d+1}{2}} \sum_\Delta \sum_\Psi 
\|\Xi_{\Psi}*(\psi_\Delta f)\|_p^p,
\]
Observe now that
\begin{equation}\label{zzz.11}
\Xi_{\Psi}*(\psi_\Delta \cdot (\Xi_{\Psi'}*f))\not\equiv 0
\end{equation}
is possible only if the Fourier supports of 
$\Xi_{\Psi}$ and $\psi_\Delta \cdot (\Xi_{\Psi'}*f))$ intersect,
i.e. if $\Gamma_{N,\Psi}(C)$ intersects the $C\sqrt{N}$-neighbourhood
of $\Gamma_{N,\Psi'}(C)$.  Since $\sqrt{N}\ll N^{3/4}$, the latter
set has (for large $N$) roughly the same size as $\Gamma_{N,\Psi'}(C)$.
Hence the number of $\Psi$ for which \bref{zzz.11} holds with
a fixed $\Psi'$ is bounded by a constant (independent of $N$ and
$\Psi'$), and similarly with $\Psi$ and $\Psi'$ interchanged.
It follows that
\begin{eqnarray*}
\sum_\Psi \sum_\Delta \|\Xi_{\Psi}\ast(\psi_\Delta f)\|_p^p
&\lesssim&\sum_{\Psi}\sum_\Delta \sum_{\Psi'}
\|\Xi_{\Psi}\ast(\psi_\Delta \cdot (\Xi_{\Psi'}\ast f))\|_p^p
\\
&\lesssim&\sum_\Delta\sum_{\Psi'}\|\psi_\Delta \cdot(\Xi_{\Psi'}\ast f)\|_p^p
\\
&\lesssim&\sum_{\Psi'}\|\Xi_{\Psi'}\ast f\|_p^p,
\end{eqnarray*}
which proves \bref{ls10}.

We now fix a $\Psi$ as above, and let $L_\Psi$ be a Lorentz
transformation mapping $\G_{N,\Psi}(C)$ to a sector of
$\G_{\sqrt{N}}(C')$.  
Namely, suppose that $\Psi$ is centered at a point
$e\in S^{d-1}$, and let $\omega_1,\dots,\omega_{d-1}$ be vectors orthogonal
to $(e, 1)$ and $(e, -1)$. Then $L_{\Psi}$ is the transformation
mapping $(e,1)$ to $N^{-1/2}(e, 1)$, $(e,-1)$ to $(e, -1)$, and
$\omega_i$ to $N^{-1/4}\omega_i$. 

Let $g_\Psi=(\Xi_\Psi * f)\circ L_\Psi$,  
then $\widehat{g_\Psi}$ is supported on $\G_{\sqrt{N}}(C')$.
We further have
\[
\|g_\Psi\|_{\infty, mic}\lesssim 1.
\]
This follows from the fact that  sectors of $\G_N$
of angular length $N^{-\frac{1}{2}}$ contained in $\Psi$ correspond
to sectors of $\G_{\sqrt{N}}$ of angular width $N^{\frac{1}{4}}$
under $L_\Psi$.

Applying the inductive hypothesis \bref{4.1} on scale
$N^{\frac{1}{2}}$ to the functions $g_\Psi$, and using
\bref{mnmnmn}, we conclude that 
\[
\|\Xi_{\Psi}\ast f\|_p^p\leq \d^{-C\e}\d^{-\frac{\alpha}{2}}
\d^{\frac{d}{2}-\frac{(d-1)p}{4}}\|\Xi_{\Psi}\ast f\|_2^2.
\]
Combining this with \bref{ls10} and using the essential orthogonality
of $\Xi_\Psi * f$, we obtain \bref{zzz.10} as claimed.
\hfill$\square$

\begin{lemma}\label{Lemma 5.2}
Assume that $P(p,\alpha)$ holds, and let $f$ be an
$N$-function with associated family of plates $\P$ which localizes
at height $\l$. Then for any $\beta>(1-\e_0)\alpha$ we have
\begin{equation}\label{4.2}
|\{ |f|>\l\}|\lesssim \l^{-p} \d^{d-\frac{(d-1)p}{2}-\beta}
\d^{\frac{d+1}{2}}|\P|.
\end{equation}
\end{lemma}

\underline{Proof} Let $W=\{|f|\geq\l\}$.  The localization
assumption means that $f$ has subfunctions $f^Q$, where $Q$
ranges over $t$-cubes, such that \bref{vo9} holds and
\[
|W|\lessapprox |\bigcup_Q W_Q|,
\]
where $ W_Q=\{x\in Q:|f_Q(x)|\gtrapprox \l\}$.

Let $g_Q=(\psi_Qf^Q)\circ a_Q$.  By Lemma \ref{Lemma 4.2} we have
$\|g_Q\|_{\infty, mic}\lesssim t^{-\frac{d-1}{2}}$.  Applying 
the inductive hypothesis \bref{4.1} to $t^{\frac{d-1}{2}}g_Q$,
with $N$ replaced by $tN$ and $\l$ replaced by $\lo{-C}
t^{\frac{d-1}{2}}\l$, we obtain 
\begin{eqnarray*}
|\{|g_Q|\gtrapprox \l\}|
&\lessapprox &\l^{-p}(tN)^{\frac{(d-1)p}{2}-d+\alpha}
t^{-(p-2)\frac{d-1}{2}}\|g_Q\|_2^2
\\
&=&\l^{-p}\d^{d-\frac{(d-1)p}{2}-\alpha}t^{-1}\|g_Q\|_2^2.
\end{eqnarray*}
By \bref{Schur}, \bref{ls3}, and \bref{vo9}, we have
\[
\sum_Q \|g_Q\|_2^2\lesssim t^{-d}\sum_Q \|f^Q\|_2^2
\lesssim \d^{\frac{d+1}{2}}\sum_Q |\P_Q|.
\lessapprox \d^{\frac{d+1}{2}}|\P|.
\]
Hence
\begin{eqnarray*}
|W|\lessapprox \sum_Q|W_Q|
&\lessapprox & t^{d+1}\sum_Q|\{|g_Q|\gtrapprox \l\}|\\
&\lessapprox &\l^{-p}\d^{d-\frac{(d-1)p}{2}-\alpha}
t^{\alpha+d}\sum_Q\|g_Q\|_2^2\\
&\lessapprox &\l^{-p}\d^{d-\frac{(d-1)p}{2}-\alpha}
t^{\alpha}\cdot\d^{\frac{d+1}{2}}|\P|.
\end{eqnarray*}
The lemma follows since $t=\d^{\e_0}$.
\hfill$\square$

\begin{lemma}\label{Lemma 5.3}
Assume that we know $P(p,\alpha)$, and let $f$ be an
$N$-function satisfying \bref{vo5}. Then \bref{4.2} holds
if $\beta>(1-\frac{\e_0}{2})\alpha$.
\end{lemma}

\underline{Proof}
Let $W=\{|f|\geq\l\}$.  By \bref{vo5}, we may apply Lemma \ref{Lemma 3.2}
to $f$ and $\l$.  If $f$ localizes at $\l$, then the conclusion follows
from Lemma \ref{Lemma 5.2}.  Otherwise, there is a subfunction $f^*$ of $f$ 
such that $|f^*|\gtrapprox \l$ on a logarithmic fraction 
$W^*$ of $W$ and that
\begin{equation}
\|\psi_{\Delta}f^*\|_2^2\leq t^{5d} \d^{\frac{3}{4}(d+1)}\l^2
\label{vx9}\end{equation}
for each $\d^{\frac{1}{2}}$-cube $\Delta$. 
We apply Lemma \ref{Lemma 5.1} to $f^*$, obtaining a family of
$\sqrt{N}$-functions $f_{\Delta}$ and a value of $\lambda_*$ as in
\bref{vp.lambda} so that
\begin{equation}\label{zzz.19}
\{|f^*|\gtrapprox \l\}|\lessapprox |\bigcup_{\Delta}
a_{\Delta}^{-1}(\{|f_{\Delta}|\geq\l_*\})|
\end{equation}
\begin{equation}\label{zzz.20}
|\P(f_{\Delta})|\leq\d^{-C\e}t^{4d}\l_*^2,
\end{equation}
\begin{equation}\label{zzz.21}
\sum_{\Delta}|\P(f_{\Delta})|\leq\d^{-C\e}\d^{-\frac{3}{4}(d+1)}
(\frac{\l_*}{\l})^{p}
\d^{\frac{d}{2}-\frac{(d-1)p}{4}-\frac{\alpha}{2}}\|f^*\|_2^2
\end{equation}
(here \bref{zzz.20} follows by substituting \bref{vx9} in
\bref{vp2}).
By \bref{zzz.20}, $f_\Delta$ satisfy \bref{vo8} with $\l_*$ and $\sqrt{N}$
replacing $\l$ and $N$.  Hence $f_\Delta$ localize at $\l_*$, and by Lemma
\ref{Lemma 5.2}
\begin{equation}\label{wb1}
|\{|f_{\Delta}|\geq\l_*\}|
\leq\d^{-C\e}\l_*^{-p} \d^{\frac{d+1}{4}}
\d^{\frac{d}{2}-\frac{(d-1)p}{4}-\frac{(1-\e_0)\alpha}{2}}
|\P(f_{\Delta})|.
\end{equation}
Substituting this in \bref{zzz.19}, we obtain
\begin{eqnarray*}
|\{|f|\geq\l\}|&\lessapprox&
\d^{\frac{d+1}{2}}\sum_\Delta |\{|f_{\Delta}|\geq\l_*\}|
\\
&\leq&\d^{-C\e}\sum_{\Delta}
\l_*^{-p}\d^{\frac{3d+3}{4}}
\d^{\frac{d}{2}-\frac{(d-1)p}{4}-\frac{(1-\e_0)\alpha}{2}}
\sum_\Delta |\P(f_{\Delta})|
\\
&\leq&\d^{-C\e}\l^{-p}
\d^{d-\frac{(d-1)p}{2}-\alpha(1-\frac{\e_0}{2})}
\|f\|_2^2.
\end{eqnarray*}
The conclusion follows since $\|f\|_2^2\lesssim \d^{\frac{d+1}{2}} |\P|$
by Lemma 4.1 with $p=2$. 
\hfill$\square$

\bigskip

\underline{Proof of Proposition \ref{scales}}
Let $f$ be a function such that $\supp\widehat{f}\subset \Gamma_N(1)$ and 
$\|f\|_{\infty, mic}\leq 1$.
We have already observed that \bref{4.2} follows from Tchebyshev's inequality if
\bref{cheb} holds.  Therefore we may assume that \bref{cheb} fails, i.e.
\begin{equation}
\l\geq\d^{-\frac{d-1}{2}+\frac{1}{p-2}}
\label{4.3}\end{equation}

Let $f_\Delta$ be the $\sqrt{N}$-functions constructed in Lemma \ref{Lemma 5.1}.
We claim that if \bref{4.3} holds and $p>p_d$, then $f_\Delta$ satisfy the
assumptions of either Lemma \ref{Lemma 5.2} or Lemma \ref{Lemma 5.3}
with $\d$ and $\l$ replaced by $\sqrt{\d}$ and $\l_*$.
Indeed, we have $\|\psi_{\Delta}f\|_2^2\lesssim\d$ by \bref{ls30} and rescaling;
plugging this into \bref{vp2} we obtain that
\begin{equation}\label{yyy.1}
|\P(f_{\Delta})|\lesssim \d^{-C\e}\frac{\l_*^2}{\l^2}\d^{-\frac{3d-1}{4}}.
\end{equation}

Assume first that $p>2+\frac{8}{d-3}$.  Then by \bref{4.3}
\[
|\P(f_{\Delta})|\lesssim \l_*^2 \d^{d-1-\frac{2}{p-2}-\frac{3d-1}{4}}
\lesssim t^C\l_*^2,
\]
where we used that $t=\d^{\e_0}$. Thus \bref{vo8} holds, and by Lemma \ref{Lemma 5.2}
$f_\Delta$ localize.

If on the other hand $p>2+\frac{32}{3d-7}$, then we have from \bref{4.3}
and \bref{vp.lambda} (after some algebra)
\[
\l^2\l_*^2\d^{\frac{9d-5}{8}}\gtrsim 
\d^{C\e}\l^4\d^{\frac{d-1}{2}+\frac{9d-5}{8}}\gtrsim t^{-C}.
\]
It follows that $\l^{-2}\lesssim t^C\l_*^2\d^{\frac{9d-5}{8}}$.  Substituting
this in \bref{yyy.1} yields
\[
|\P(f_{\Delta})|\lesssim \l_*^4 t^C\d^{\frac{3d-3}{8}},
\]
which is \bref{vo5} on scale $\sqrt{\d}$.

We may thus apply Lemma \ref{Lemma 5.2} or \ref{Lemma 5.3} respectively
to $f_{\Delta}$, and obtain that
\[
|\{|f_{\Delta}|\geq\l^*\}|\lesssim 
\l_*^{-p} \d^{\frac{d+1}{4}}
\d^{\frac{d}{2}-\frac{(d-1)p}{4}-\frac{\gamma}{2}} |\P(f_{\Delta})|.
\]
for any $\gamma>(1-\frac{\e_0}{2})\alpha$. Hence
\begin{eqnarray*}
|\{|f|\geq\l\}|&\lessapprox &
\d^{\frac{d+1}{2}}\sum_\Delta |\{|f_{\Delta}|\geq\l_*\}|
\\
&\lesssim &\d^{-C\e}\d^{\frac{d+1}{2}}\sum_{\Delta}
\l_*^{-p} \d^{\frac{d+1}{4}}
\d^{\frac{d}{2}-\frac{(d-1)p}{4}-\frac{\gamma}{2}} |\P(f_{\Delta})|
\\
&\leq&\d^{-C\e}
\l^{-p}\d^{d-\frac{(d-1)p}{2}-\frac{\alpha+\gamma}{2}} \|f\|_2^2
\end{eqnarray*}
as required.
\qed


\noindent{\sc Department of Mathematics, University of British Columbia, Vancouver, BC, 
Canada V6T 1Z2}

\noindent{\it ilaba@math.ubc.ca}

\bigskip

\noindent{\sc Department of Mathematics, California Institute of Technology, Pasadena, 
CA 91125, USA}


\end{document}